\documentclass[a4paper,12pt]{article}
\usepackage[T2A]{fontenc}
\usepackage{amsfonts,amssymb,amscd,amsmath}
\usepackage[dvips]{graphicx}
\usepackage[right=25mm,left=25mm,top=30mm,
bottom=40mm]{geometry}
\usepackage[cp1251]{inputenc}
\begin{document}

\begin{center}{\bf Kinks in the
parametrically excited sine-Gordon 
equation and method of averaging}\end{center}

\vskip 14pt\noindent

\centerline{Vladimir~Burd}

\vskip 14pt\noindent

Department of Mathematics, Yaroslavl State
University, Russia

\vskip 14pt\noindent

Abstract 

Parametrically excited sine-Gordon
equation is considered.
Excitation is a fast oscillating
periodic function with zero mean.
Technique of classical method of
averaging enables to 
construct the averaged equations 
in a variety of assumptions 
about driving amplitude. 
The averaged equation possesses
kinks solutions.
The results can be applied to 
the study of movement of Bloch walls 
for ferromagnetic crystals in 
the presence of a rapidly 
oscillating
magnetic field and to describe
the fluxon dynamics in long
Josephson junctions.

Keywords

sine-Gordon
equation, parametric 
excitation,
method of averaging, kinks.

%\MSC[2010] 34A37 \sep 34C29

\section{Introduction}
\label{s}
As is well known, the 
sine-Gordon equation has a 
different physical applications.
The sine-Gordon equation arose 
in the study of wave  in 
quasi-one-dimensional  
ferromagnetic materials
and the propagation 
of spin waves in superfluid 
phases $A$ and $B$ of helium 
$^3He$.
The sine-Gordon equation has 
been studied 
in connection with the 
Josephson junctions in the 
theory of superconductivity.
\par
A significant number of studies 
have been devoted to development 
of technique of perturbation 
theory for systems close 
to integrable, in particular 
to systems close to the equation
sine-Gordon [1--6]. 
The dynamics of kinks under the
action of parametric 
perturbation was considered in
[6,7,8]. 
%The action
%of parametric excitation on
%sine-Gordon equation
% was considered in [6,7,8]. 
This equation has a form
$$
u_{tt}(x,t)-u_{xx}(t,x)+f\left(\frac{t}{\varepsilon}\right)
\sin u(t,x)=0. \eqno(1)
$$
Here $f(t)$ is a periodic function
with zero mean and a constant
amplitude, $\varepsilon$ is
small positive parameter.
In the papers [7, 8] was found 
averaged equation. Averaged equation
possesses $\pi$-kinks solutions.
%Apparently results can be 
%applied to quasi-homogeneous 
%ferromagnets, which are under 
%the acting of rapidly 
%oscillating external magnetic 
%field. The presence of kinks 
%affects the motion of the walls 
%of ferromagnet.
\par
The problem of constructing the 
averaged equation is solved 
in [7,8] as 
follows. It is introduced the
Hamiltonian
$$
H=\int\limits_{-\infty}^\infty
\left(\frac{p^2}{2}
+\frac{u_x^2}{2}-f\left(\frac{t}{\varepsilon}\right)
\cos u\right)dx,
$$
where $p=u_t$. Then, 
the series of canonical near-identical
transformations is applied.
It allows to remove the fast
oscillating terms of lower orders
from the Hamiltonian.
Another way of constructing the 
averaged dynamics was proposed 
in [6]. The solution 
parametrically 
excited sine-Cordon equation 
is sought in the form of a 
Fourier series with slowly 
varying coefficients on the time
scale $\omega^{-1}$, $\omega$
being the frequency of the rapidly
varying perturbation. .
\par
In this paper we propose a 
third method of constructing of
averaged equation for (1).
It is the classical method 
of averaging (see [9,10]).
Our approach is similar to 
the method that was used in the
study of the stability of the
upper equilibrium of the pendulum 
with a vertically vibrating pivot
[11].

\section{Constructing of 
averaged equations}
\label{t}
We denote by $f_{-1}(t)$ the
periodic function with zero mean 
value whose derivative 
$f_{-1}'(t)$
satisfies the equality
$$
f_{-1}'(t)=f(t).
$$
From equation (1) we move 
to equivalent system 
$$
\begin{array}{l}
u_t=p-\varepsilon f_{-1}\left(
\frac{t}{\varepsilon}\right)
\sin u\\
p_t=u_{xx}+\varepsilon f_{-1}\left(
\frac{t}{\varepsilon}\right)p
\cos u-
\frac{1}{2}\varepsilon^2
\left(f_{-1}\left(
\frac{t}{\varepsilon}\right)
\right)^2\sin 2u.
\end{array}\eqno(2)
$$
This transition is essentially 
a transition from a recording 
of system in the form of Lagrange 
to recording of system in the 
Hamiltonian form. In system
(2) turn to fast 
time $\tau=\frac{t}
{\varepsilon}$. We 
obtain a system
$$
\begin{array}{l}
u_\tau=\varepsilon p-
\varepsilon^2f_{-1}(\tau)
\sin u\\
p_\tau=\varepsilon u_{xx}+
\varepsilon^2f_{-1}(\tau)
p\cos u-\frac{1}{2}\varepsilon^3
(f_{-1}(\tau))
^2\sin 2u.
\end{array}\eqno(3)
$$
The right-hand sides of the
system (3) are
proportional to small parameter
$\varepsilon$.
This is a standard form 
for applying the method 
of averaging. We make
regular substitution of
method of averaging 
$$
\begin{array}{l}
u=\xi+\varepsilon^2
v_2(\tau,\xi)+
\varepsilon^3v_3(\tau,\xi),\\
p=\eta+\varepsilon^2w_2(\tau,\xi,\eta)+
\varepsilon^3w_3(\tau,\xi,\eta).
\end{array}\eqno(4)
$$
The change should 
eliminate a variable 
$\tau$ from the right-hand 
side of the system (3) up to
terms of the fourth order. 
%of smallness. 
Therefore we need to get the 
system
$$
\begin{array}{l}
\xi_\tau=\varepsilon\eta+\varepsilon^2A_2(\xi,\eta)+
\varepsilon^3A_3(\xi,\eta)+O(\varepsilon^4), 
\\
\eta_\tau=\varepsilon \xi_{xx}+\varepsilon^2B_2(\xi,
\eta)+\varepsilon^3B_3(\xi,\eta)+O(\varepsilon^4).
\end{array}\eqno(5).
$$
Let us find coefficients
$A_2$,~$A_3$,~$B_2$,~$B_3$
of system (5). By 
substituting (4)
into (3) and 
replacing $\xi_{\tau}$ and
$\eta_{\tau}$ with
right-hand sides of the
system (5) we get
$$
\begin{array}{l}
\varepsilon\eta+\varepsilon^2
A_2(\xi,\eta)+
\varepsilon^3A_3(\xi,\eta)+
\varepsilon^2v_{2\tau}(\tau,\xi)
+\varepsilon^3v_{2\xi}(\tau,\xi)
\eta
+\varepsilon^3
v_{3\tau}(\tau,\xi)+O(\varepsilon^4)=\\
\varepsilon\eta-\varepsilon^2
f_{-1}(\tau)
\sin \xi+\varepsilon^3w_2(\tau,
\xi,\eta)+O(\varepsilon^4),
\\
\varepsilon \xi_{xx}+
\varepsilon^2B_2(\xi,
\eta)+\varepsilon^3B_3(\xi,\eta)+
\varepsilon^2w_{2\tau}(\tau,\xi,\eta)
+\\
\varepsilon^3
w_{2\xi}(\tau,\xi,\eta)\eta
+\varepsilon^3w_{2\eta}(\tau,\xi,\eta)
\xi_{xx}+\varepsilon^3
w_{3\tau}(\tau,\xi,\eta)+
O(\varepsilon^4)\\
=\varepsilon\xi_{xx}+\varepsilon^3
\frac{\partial^2}{\partial x^2}
(v_2(\tau,\xi))
+\varepsilon^2f_{-1}(\tau)\eta
\cos\xi-\varepsilon^3
\frac{1}{2}(f_{-1}(\tau))^2
\sin 2\xi+O(\varepsilon^4).
\end{array}\eqno(6)
$$
Equating the coefficients
of power of $\varepsilon^2$
we obtain
$$
v_{2\tau}(\tau,\xi)+A_2(\xi,\eta)=-f_{-1}
(\tau)\sin \xi,\quad
w_{2\tau}(\tau,\xi,\eta)+B_2(\xi,\eta)
=f_{-1}(\tau)\eta\cos\xi.
$$
The function $A_2(\xi,\eta)$
is defined as the mean 
value on $\tau$ of the 
function $-f_{-1}(\tau)
\sin \xi$. 
Hence we obtain
$A_2(\xi,\eta)\equiv 0$.
Then the function 
$v_2(\tau,\xi)$ is uniquely
determined as periodic 
function on $\tau$
with zero mean value.
The function is defined by formula
$v_2(\tau,\xi)=-
(\int f_{-1}(s)ds)
\sin \xi$. Similarly we find
that $B_2(\xi,\eta)\equiv 0$. 
The function 
$w_2(\tau,\xi,\eta)$ 
is periodic function on 
$\tau$ with zero mean 
value and
$w_2(\tau,\xi,\eta)=
(\int f_{-1}(s)ds)\eta
\cos\xi$.
Equating coefficients of
power $\varepsilon^3$ in
(6) we obtain
$$
A_3(\tau,\xi)+v_{3\tau}(\tau,\xi)
+v_{2\xi}(\tau,\xi)
\eta+v_{2\eta}(\tau,\xi)\xi_{xx}
=w_2(\tau,\xi,\eta),
$$
$$
B_3(\xi,\eta)+w_{3\tau}(\tau,\xi,\eta)+
w_{2\xi}(\tau,\xi,\eta)
\eta+w_{2\eta}(\tau,\xi,\eta)
\xi_{xx}=
-\frac{1}{2}(f_{-1}
(\tau))^2\sin 2\xi+
\frac{\partial^2}
{\partial x^2}(v_2(\tau,
\xi)).
$$
%Then $A_3(\xi,\eta)\equiv 0$ as
Mean values of functions 
$v_{2\xi}(\tau,\xi)\eta$,~
$v_{2\eta}(\tau,\xi)\xi_{xx}$,~
$w_2(\tau,\xi,\eta)$ 
are zero. Then $A_3(\xi,\eta)
\equiv 0$.
The 
same arguments shows that
$B_3(\tau,\xi,\eta)$ is 
determined by formula 
$$
B_3(\xi,\eta)=\frac{1}{2}
\langle(f_{-1}(\tau))^2\rangle
\sin 2\xi,
$$
where $\langle f_{-1}^2
(\tau)\rangle$ is mean
value of periodic function
$f_{-1}^2(\tau)$.
Therefore the averaged
system has the form
$$
\begin{array}{l}
\xi_\tau=\varepsilon\eta,\\
\eta_\tau=\varepsilon
\xi_{xx}-
\varepsilon^3
\frac{1}{2}\langle(f_{-1}(\tau))^2\rangle
\sin 2\xi.
\end{array}
$$
Going back to the time
$t=\varepsilon\tau$ we
obtain the system
$$
\begin{array}{l}
\xi_t=\eta,\\
\eta_t=\xi_{xx}-
\varepsilon^2\frac{\Delta}{2}
\sin2\xi,
\end{array}
$$
where
$$
\Delta=\langle(f_{-1}(\tau))^2\rangle.
$$
This system can be 
written as an equation 
of the second-order
$$
\xi_{tt}-\xi_{xx}+\varepsilon^2
\frac{\Delta}{2}
\sin 2\xi=0. \eqno(7)
$$
The equation (7) has the
solutions
$$
u(x,t)=2\arctan\left
[\exp\left((x-ct)\frac
{\varepsilon\sqrt
{\Delta}}{\sqrt{1-c^2}}+
\delta\right)\right],
$$
where $\delta$ is constant.
These solutions are $\pi$
kinks. 
%We believe that
%kinks are approximate 
%solutions of equation (1) 
%for sufficiently small
%$\varepsilon$.
\par
It was assumed above 
that the function 
$f(t/\varepsilon)$
has a constant amplitude.
Now we assume that amplitude of
excitation have the order of
$1/\varepsilon$. Let
excitation is 
$$
\frac{1}{\varepsilon}f\left
(\frac{t}{\varepsilon}\right).
$$
Parametric excitation of
the sine-Gordon equation
now has the form
$$
u_{tt}(x,t)-u_{xx}
(t,x)+\frac{1}
{\varepsilon}
f\left(\frac{t}
{\varepsilon}\right)
\sin u(t,x)=0. \eqno(8)
$$
We again denote by $f_{-1}(t)$ 
the periodic function with
zero mean value that derivative
$f'_{-1}(t)$ satisfies the 
equality
$$
f_{-1}'(t)=f(t).
$$
From equation (8) we 
pass to the equivalent 
system of equations
$$
\begin{array}{l}
u_t=p-f_{-1}\left(
\frac{t}{\varepsilon}
\right)\sin u\\
p_t=u_{xx}+f_{-1}\left(
\frac{t}{\varepsilon}\right)p
\cos u-
\frac{1}{2}\left(f_{-1}
\left(\frac{t}{\varepsilon}\right)
\right)^2\sin 2u.
\end{array}
$$
Then move on to the 
fast time $\tau=t/\varepsilon$. 
We obtain the system
$$
\begin{array}{l}
u_\tau=\varepsilon(p-
f_{-1}(\tau)
\sin u)\\
p_\tau=\varepsilon
(u_{xx}+f_{-1}(\tau)
p\cos u-\frac{1}{2}
(f_{-1}(\tau))^2
\sin 2u).
\end{array}
$$
The averaged system is 
defined in a first approximation
and has the form
$$
\begin{array}{l}
\xi_\tau=\varepsilon \eta\\
\eta_\tau=\varepsilon(\xi_{xx}-\frac{1}{2}
\langle(f_{-1}(\tau))^2\rangle\sin 2\xi).
\end{array}
$$
We assume $\Delta=
\langle f_1(\tau)^2\rangle$ and 
record the average system 
in the form of 
a second-order equation in
time $t$
$$
\xi_{tt}-\xi_{xx}+
\frac{\Delta}{2}\sin 
2\xi=0. \eqno(9).
$$
The equation (9) has the $\pi$-
kinks
$$
u(x,t)=2\arctan\left[\exp\left((x-ct)\frac
{\sqrt{\Delta}}{\sqrt{1-c^2}}+\delta\right)
\right],\quad  \delta=const.
$$
\par
Let's look at the 
excitement of the 
following form
$$
u_{tt}(x,t)-u_{xx}(t,x)+
\left(1+\frac{1}{\varepsilon}
f\left(\frac{t}{\varepsilon}\right)\right)\sin u(t,x)
=0. \eqno(10)
$$
Note that in the case of 
one variable, equation (10)
is transformed into the 
equation of motion of 
the pendulum with a 
vertically oscillating 
pivot and 
small amplitude. In the 
notations introduced earlier 
system of equations is 
equivalent to the 
equation (10) has the form
$$
\begin{array}{l}
u_t=p-f_{-1}\left(
\frac{t}{\varepsilon}\right)\sin u,\\
p_t=u_{xx}-\sin u+f_{-1}\left(
\frac{t}{\varepsilon}\right)p\cos u    
-\frac{1}{2}\left(f_{-1}\left(
\frac{t}{\varepsilon}
\right)\right)^2\sin 2u.
\end{array}\eqno(11)
$$
We go to a fast time 
$\tau=t/\varepsilon$. 
Then we average system 
(11). Averaged system 
of the first approximation 
in time $t$ can be 
written as equation of second order
$$
\xi_{tt}-\xi_{xx}+
\sin\xi+\frac{\Delta}{2}
\sin 2\xi=0. \eqno(12)  
$$
where $\Delta=
\langle f_1(\tau)^2\rangle$.
Thus averaged equation 
is the double sine-Gordon
equation. The kinks solution
of equation (12) may be
written in the form
$$
\xi(t,x)=2\arctan\left[
\frac{1}{\sqrt{1+\Delta/2}}
csch\left(\sqrt{1+
\Delta/2}\frac{x-ct}
{\sqrt{1-c^2}}\right)\right].
$$
\par
Similarly, the following 
parametric perturbed equation
$$
u_{tt}-u_{xx} +\left
(1+f\left(
\frac{t}{\varepsilon}
\right)\right)
\sin u=0.   \eqno(13)
$$
is considered. The equivalent 
system of equations is
$$
\begin{array}{l}
u_t=p-\varepsilon f_{-1}\left(
\frac{t}{\varepsilon}\right)
\sin u\\
p_t=u_{xx}-\sin u+\varepsilon f_{-1}\left(
\frac{t}{\varepsilon}\right)p
\cos u-
\frac{1}{2}\varepsilon^2
\left(f_{-1}\left(
\frac{t}{\varepsilon}\right)
\right)^2\sin 2u.
\end{array}
$$
The averaged equation for (13)
has the form
$$
\xi_{xx}-\xi_{xx}+\sin\xi+
\varepsilon^2\frac{\Delta}{2}
\sin 2\xi=0. 
$$

\section{Conclusion}
\label{u}
The classical method of averaging 
is proposed for the 
construction of the averaged 
dynamics in parametrically 
perturbed sine-Gordon equation. 
The essential point is transition
to an equivalent system of 
first order equations. The 
averaging method is also 
applicable to the case of the 
direct driving force in the 
sine-Gordon equation and driven 
damped sine-Gordon model.
This method can be used to 
analyze the dynamics of kinks
on rotationg and oscillating 
background.

\end{document}